\documentclass[12pt,leqno]{article}

\usepackage{amsfonts}
\usepackage{amssymb}

\newcommand{\cE}{{\cal E}}
\newcommand{\cF}{{\cal F}}

\newcommand{\cK}{{\cal K}}
\newcommand{\cL}{{\cal L}}

\newcommand{\cO}{{\cal O}}

\newcommand{\ord}{{\rm ord}}
\newcommand{\End}{{\rm End}}
\newcommand{\Hom}{{\rm Hom}}
\newcommand{\Aut}{{\rm Aut}}
\newcommand{\id}{{\rm id}}

\newcommand{\Ind}{{\rm Ind}}
\newcommand{\Res}{{\rm Res}}

\newcommand{\Trace}{{\rm Trace}}

\newcommand{\ZZ}{{\mathbb Z}}
\newcommand{\CC}{{\mathbb C}}
\newcommand{\FF}{{\mathbb F}}

\newcommand{\NN}{{\mathbb N}}

\newcommand{\QQ}{{\mathbb Q}}

\newcommand{\ra}{\rightarrow}
\def\rightepi{{\longrightarrow \kern-0.7em \rightarrow}}
\newcommand{\oplusm}{\mathop{\oplus}\limits}
\newcommand{\summ}{\mathop{\sum}\limits}

\newcommand{\notteilt}{{\,\not{\kern-0.075em|}\,}}
\def\antiddots{\mathinner{\mkern1mu\raise1pt\vbox{\kern7pt\hbox{.}}\mkern2mu
    \raise4pt\hbox{.}\mkern2mu\raise7pt\hbox{.}\mkern1mu}}

\begin{document}

\vspace*{15ex}

\begin{center}
{\LARGE\bf Computing the equivariant Euler characteristic of
Zariski and \'etale sheaves on curves}\\
\bigskip
by\\
\bigskip
{\sc Bernhard K\"ock}\\
\bigskip
{\em Dedicated to Vic Snaith on the occasion of his 60th birthday}
\end{center}

\bigskip

\begin{quote}
{\footnotesize {\bf Abstract}. We prove an equivariant
Grothendieck-Ogg-Shafarevich formula. This formula may be viewed
as an \'etale analogue of well-known formulas for Zariski sheaves
generalizing the classical Chevalley-Weil formula. We give a new
approach to those formulas (first proved by Ellingsrud/L{\o}nsted,
Nakajima, Kani and Ksir) which can also be applied in the \'etale
case.


}

\end{quote}

\bigskip

\section*{Introduction}

This paper deals with the Riemann-Roch problem for equivariant
Zariski sheaves and equivariant \'etale sheaves on smooth
projective curves, i.e.\ with the computation of their
(equivariant!) Euler characteristic. In the case of Zariski
sheaves we give a new, very natural and quick approach to
generalizations of the classical Chevalley-Weil formula proved by
Ellingsrud/L{\o}nsted (see \cite{EL}), Kani (see \cite{Ka}) and
Nakajima (see \cite{Na}) in the 1980s and we derive
generalizations of a more recent result of Ksir (see \cite{Ks}).
In the case of \'etale sheaves we develop and prove an equivariant
Grothendieck-Ogg-Shafarevich formula by imitating our new approach
for Zariski sheaves.

Let $X$  be a connected smooth projective curve over an
algebraically closed field $k$ and let $G$ be a finite subgroup of
$\Aut(X/k)$ of order $n$. We assume throughout this paper that the
canonical projection $\pi: X \ra Y:=X/G$ is tamely ramified.

Using the coherent Lefschetz fixed point formula (see \cite{Do1}
or \cite{BFQ} or \cite{Ko}) in conjunction with the Riemann-Roch
formula and Hurwitz formula, we prove the following theorem (see
Theorem 1.1) which effectively implies all known formulas (see
Corollaries 1.3, 1.4, 1.7 and 1.8) for the equivariant Euler
characteristic
\[\chi(G,X,\cE):=[H^0(X,\cE)] - [H^1(X,\cE)] \]
of a locally free $G$-sheaf $\cE$ on $X$, considered as an element
of the Grothendieck group $K_0(G,k)$ of all $k$-representations of
$G$.

{\bf Theorem}.
\[\chi(G,X,\cE) =  \left((1-g_Y)r + \frac{1}{n}\deg(\cE)
\right)[k[G]] - \frac{1}{n}\sum_{P\in X} \sum_{d=0}^{e_P-1}
d[\Ind_{G_P}^G(\cE(P)\otimes \chi_P^d)].\]

Here, $r$ denotes the rank of $\cE$, $\cE(P)$ the fibre of $\cE$
at $P$ and $\chi_P$ the character of $G_P$ which is given by the
representation of $G_P$ on the cotangent space ${\mathfrak
m}_P/{\mathfrak m}_P^2$.

We now pass to the \'etale world. We fix a prime $l \not= {\rm
char}(k)$. Let $\cF$ be a constructible $\FF_l$-sheaf on the
\'etale site $X_{\rm \acute{e}t}$ which carries a $G$-action
compatible with the given $G$-action on $X$. We are interested in
computing the equivariant Euler characteristic
\[\chi_{\rm \acute{e}t}(G,X,\cF) := [H^0_{\rm \acute{e}t}(X,\cF)] -
[H^1_{\rm \acute{e}t}(X,\cF)] + [H^2_{\rm \acute{e}t}(X,\cF)] \]
considered as an element of the Grothendieck group $K_0(G,\FF_l)$
of $\FF_l$-representations of $G$.

In the extreme case that $G$ is the trivial group, this problem is
solved by the classical Grothendieck-Ogg-Shafarevich formula (see
Theorem 2.12 on p.\ 190 in \cite{Mi}). In the extreme case that
$\cF$ is the constant sheaf with trivial $G$-action, a
satisfactory answer to this problem follows from Remark 2.9 on p.\
187 in \cite{Mi} (see Remark 2.2(b)). By imitating our approach
for Zariski sheaves we prove the following result for an arbitrary
group $G$ and arbitrary sheaf~$\cF$ (see Theorem 2.1):

{\bf Theorem}. {\em We assume that the characteristic of $k$ does
not divide $n$. Then we have:
\begin{eqnarray*}
\lefteqn{\chi_{\rm \acute{e}t}(G,X,\cF) =}\\
 & = &\left((2-2g_Y)\dim_{\FF_l}(\cF_{\bar{\eta}}) - \frac{1}{n}
 \alpha(\cF) + \sum_{P\in X} \frac{e_p}{n}(\dim_{\FF_l}(\cF_P) -
 \dim_{\FF_l}(\cF_{\bar{\eta}}))\right)[\FF_l[G]] \\
 && - \sum_{P\in X} \frac{e_P}{n}[\Ind_{G_P}^G(\cF_P \otimes
 I_P)].
\end{eqnarray*}}

Here, $g_Y$ denotes the genus of $Y:= X/G$, $\eta$ the generic
point of $X$, $\alpha(\cF):= \sum_{P\in X} \alpha_P(\cF)$ the sum
of the wild conductors of $\cF$, $e_P$ the ramification index of
the canonical projection $\pi: X \ra Y$ at $P \in X$ and $I_P :=
{\rm ker}(\FF_l[G_P] \ra \FF_l)$ the augmentation representation
of the decomposition group $G_P$ at $P$.

As a corollary, we obtain that $\alpha(\cF)$ is divisible by $n$
(see Corollary 2.4). Furthermore, this formula has the following
simple shape, if $\pi$ is \'etale (see Remark 2.2(a)):
\[\chi_{\rm \acute{e}t}(G,X,\cF) = \frac{1}{n}\chi_{\rm \acute{e}t}(X,\cF)[\FF_l[G]].\]
In fact, this formula is valid without the assumption that ${\rm
char}(k)$ does not divide $n$. In particular, the
(non-equivariant) Euler characteristic $\chi_{\rm
\acute{e}t}(X,\cF)$ is divisible by $n$.

{\bf Acknowledgments}. I would like to thank Igor Zhukov for
raising the equivariant Riemann-Roch problem for \'etale sheaves
and for drawing my attention to the paper \cite{Do2}. Furthermore,
I would like to thank him and Victor Snaith for helpful and
encouraging discussions, for reading carefully a preliminary
version of this paper and for suggesting several corrections.

\bigskip

\section*{\S 1 The formulas of Ellingsrud/L{\o}nsted,
Nakajima, Kani and Ksir revisited}

The object of this section is to give a new approach to
generalizations of the classical Chevalley-Weil formula published
by Ellingsrud/L{\o}nsted, Nakajima and Kani and to derive
generalizations of a comparatively simple formula recently
published by Ksir.

Let $X$ be a connected smooth projective curve over an
algebraically closed field $k$ and let $G$ be a finite subgroup of
$\Aut(X/k)$ of order $n$. We assume in this section that the
canonical projection $\pi:X \ra Y := X/G$ is tamely ramified.

We denote the genus of $X$ and $Y$ by $g_X$ and $g_Y$,
respectively. For any (closed) point $P \in X$, let $G_P:=
\{\sigma \in G: \sigma(P) = P\}$ denote the decomposition group
and let $e_P$ denote the ramification index of $\pi$ at $P$. It is
well-known (see Corollaire~1 of Proposition~7, Chapitre~IV on
p.~75 of \cite{Se1}) that $G_P$ is a cyclic group of order $e_P$
and that
\[\Hom(G_P, k^\times) = \{\chi_P^0, \ldots, \chi_P^{e_P-1}\}\]
where $\chi_P: G_P \ra k^\times$ denotes the character which is
given by the action of $G_P$ on the cotangent space ${\mathfrak
m}_P/{\mathfrak m}_P^2$. For $Q \in Y$ we set $e_Q := e_P$ where
$P \in \pi^{-1}(Q)$.

We denote the Grothendieck group of all $k$-representations of $G$
(of finite dimension) by $K_0(G,k)$. It is free with basis
$\hat{G}$ where $\hat{G}$ denotes the set of isomorphism classes
of irreducible $k$-representations of $G$.

Now, let $\cE$ be a locally free $G$-module on $X$ of rank $r$,
i.e., we have $\cO_X$-isomorphisms $g^*(\cE)\ra \cE$, $g\in G$,
which satisfy the usual composition rules. Then, the Zariski
cohomology groups $H^i(X, \cE)$, $i\ge 0$, are $k$-representations
of $G$. Let
\[\chi(G,X,\cE) := [H^0(X,\cE)] - [H^1(X,\cE)] \in K_0(G,k)\]
denote the equivariant Euler characteristic of $X$ with values in
$\cE$. For any $P\in X$, we view the fibre
$\cE(P):=\cE_P/{\mathfrak m}_P\cE_P$ as a $k$-representation of
$G_P$.

The following theorem computes the equivariant Euler
characteristic $\chi(G,X,\cE)$.

{\bf Theorem 1.1}. {\em We have in $K_0(G,k)$:
\[n \cdot \chi(G,X,\cE) = \left(n(1-g_Y)r + \deg(\cE)\right)[k[G]]
- \sum_{P\in X} \sum_{d=0}^{e_P-1} d\left[\Ind_{G_P}^G(\cE(P)
\otimes \chi_P^d)\right].\]}

{\bf Proof}. By classical representation theory (see
Corollary~(17.10) on p.~424 in \cite{CR2}) it suffices to show
that the Brauer characters of both sides of formula~(1) coincide.
For any $k$-representation $V$ of $G$ and for any $\sigma \in G$
of order prime to $p$ we write $\Trace(\sigma|V)$ for the value of
the Brauer character of $V$ at $\sigma$. Recall that
\[\Trace(\sigma|V) = \sum_{i=1} ^{\dim(V)} \varphi(\alpha_i)\]
where $\alpha_i$, $i=1, \ldots, \dim(V)$, are the eigenvalues of
the $k$-linear map $\sigma$ on $V$ and $\varphi: k^\times \ra
K^\times$ is the Teichm\"uller character from the group of
invertible elements in $k$ to the group of invertible elements in
the quotient field $K:=Q(W(k))$ of the Witt ring $W(k)$ of $k$.
(We set $K:=k$ and $\varphi:= \id$, if $p=0$.)\\
Let $\sigma \in G$ such that ${\rm char}(k)$ does not divide the
order of $\sigma$ and let $X^\sigma := \{P\in X: \sigma(P)=P\}$
denote the set of points in $X$ fixed by $\sigma$. Then we have:
\begin{eqnarray*}
\lefteqn{\sum_{P\in X} \sum_{d=0}^{e_P-1} d \cdot \Trace\left(
\sigma|\Ind_{G_P}^G(\cE(P) \otimes \chi_P^d)\right) = }\\ && =
\sum_{P\in X} \sum_{d=0}^{e_P-1} \frac{d}{e_P} \sum_{\tau\in G,\,
\tau^{-1} \sigma \tau \in G_P} \Trace(\tau^{-1}\sigma\tau |
\cE(P)) \cdot \chi_P^d(\tau^{-1}\sigma\tau)\\ && \hspace*{4cm}
\textrm{(by Lemma~(21.28) on p.~509 in \cite{CR2})}\\ &&
=\sum_{P\in X} \sum_{d=0}^{e_P-1} \frac{d}{e_P}\sum_{\tau\in G,\,
\tau(P) \in X^\sigma} \Trace(\sigma | \cE(\tau(P))) \cdot
\chi_{\tau(P)}^d(\sigma) \\ && = \sum_{P\in X^\sigma}
\frac{n}{e_P} \cdot \Trace(\sigma|\cE(P)) \cdot \sum_{d=0}^{e_P-1}
d \cdot \chi_P^d(\sigma) \\ && = \left\{
\begin{array}{ll}
n \cdot r \cdot \summ_{P\in X} \frac{e_P-1}{2} & \textrm{for }
\sigma =\id
\\ n \cdot \summ_{P\in X^\sigma} \Trace(\sigma| \cE(P)) \cdot
(\chi_P(\sigma)-1)^{-1} & \textrm{for } \sigma \not=\id \textrm{
(see Lemma 1.2 below)}
\end{array}\right.
\end{eqnarray*}
For $\sigma \not= \id$ we have $\Trace(\sigma|k[G]) = 0$, so the
character value of the right hand side of the formula in
Theorem~1.1 at the place $\sigma$ equals
\[n\cdot \sum_{P\in X^{\sigma}} \Trace(\sigma|\cE(P)) \cdot
(1-\chi_P(\sigma))^{-1}.\] By the Lefschetz fixed point formula
(see Example 3 in \cite{Ko} or \cite{BFQ} or \cite{Do1}; here we
use the assumption that ${\rm char}(k)$ does not divide
$\ord(\sigma)$), this equals the character value of the left hand
side at the place $\sigma$.\\ For $\sigma =\id$ we have:
\begin{eqnarray*}
\lefteqn{(n(1-g_Y)r + \deg(\cE)) \cdot \Trace(\sigma|k[G])=} \\ &&
= n(n(1-g_Y)r + \deg(\cE)) \\&& = n\left((1-g_X)r + \deg(\cE) + r
\sum_{P\in X}\frac{e_P-1}{2}\right)
\end{eqnarray*}
by the Hurwitz formula (see Corollary 2.4 on p.\ 301 in
\cite{Ha}). Hence, the character value of the right hand side of
the formula in Theorem~1.1 at the place $\sigma = \id$ equals
\[n((1-g_X)r + \deg(\cE)).\]
By the Riemann-Roch formula (see \S 1 in Chapter IV of \cite{Ha}
and Exercise 6.11 on p.\ 149 in \cite{Ha}), this equals the
character value of the left hand side at the place $\sigma =
\id$.\\ Thus, the proof of Theorem~1.1 is complete. \hfill $\Box$

{\bf Lemma 1.2}. {\em Let $m\in \NN$ and $\zeta \not= 1$ an $m$-th
root of unity. Then we have:
\[ m (\zeta-1)^{-1} = \sum_{d=1}^{m-1} d\zeta^{d}.\]}

{\bf Proof.} $(\sum_{d=1}^{m-1} d\zeta^d)(\zeta -1) =
\sum_{d=1}^{m-1} d\zeta^{d+1} - \sum_{d=1}^{m-1}d\zeta^d =
(m-1)\zeta^m - \sum_{d=1}^{m-1}\zeta^d =m$. \hspace*{\fill} $\Box$

{\em Remark.} A generalization of Theorem~1.1 and the subsequent
Corollary~1.4 to the so-called weakly ramified case can be found
in \cite{Ko2}.

The following corollary is the main result of the paper \cite{EL}
by Ellingsrud and L{\o}nsted; it computes the multiplicity of any
irreducible representation $V\in \hat{G}$ in the equivariant Euler
characteristic $\chi(G,X,\cE)$, if ${\rm char}(k)$ does not divide
$n$. While the proof of Ellingsrud and L{\o}nsted is based on the
study of the cokernel of the natural embedding
$\pi^*(\pi_*^G(\cE)) \hookrightarrow \cE$, we derive it from
Theorem 1.1 and hence from the Lefschetz fixed point formula.

{\bf Corollary 1.3} (Formula (3.7) in \cite{EL}). {\em We assume
that ${\rm char}(k)$ does not divide $n$. For $Q\in Y$, $d\in \NN$
and $V\in \hat{G}$, let $n_{d,Q,\cE,V}$ denote the multiplicity of
$\chi_P^d$ in $\cE(P) \otimes \Res_{G_P}^G(V)$ where $P\in
\pi^{-1}(Q)$. Then we have in $K_0(G,k)_\QQ$:
\[\chi(G,X,\cE) = \left(\frac{1}{n} \deg(\cE) +
(1-g_Y)r\right)[k[G]] - \sum_{V\in \hat{G}} \sum_{Q\in Y}
\sum_{d=1}^{e_Q}\left(1-\frac{d}{e_Q}\right)n_{d,Q,\cE,\check{V}}[V].\]}

{\bf Proof.} Let $\langle\, ,\, \rangle: K_0(G,k) \times K_0(G,k)
\ra \ZZ$ denote the usual character pairing. Then, for all $P\in
X$ and $d\in \{0,\ldots,e_P-1\}$, we have:
\begin{eqnarray*}
\lefteqn{\sum_{V\in \hat{G}}\langle[\chi_P^d], [\cE(P) \otimes
\Res_{G_P}^G(\check{V})]\rangle[V] = }\\&& = \sum_{V\in \hat{G}}
\langle\langle[\Ind_{G_P}^G(\cE(P) \otimes
\chi_P^{-d})],[V]\rangle[V]\\&& = [\Ind_{G_P}^G(\cE(P) \otimes
\chi_P^{-d})] \quad \textrm{in} \quad K_0(G,k).
\end{eqnarray*}
Hence we have:
\begin{eqnarray*}
\lefteqn{\sum_{V\in\hat{G}}\sum_{Q\in
Y}\sum_{d=1}^{e_Q}\left(1-\frac{d}{e_Q}\right)
n_{d,Q,\cE,\check{V}} [V] =} \\ && = \sum_{P\in X} \frac{e_P}{n}
\sum_{d=1}^{e_P-1} \left(1-\frac{d}{e_P}\right)\sum_{V\in\hat{G}}
\langle[\chi_P^d],[\cE(P)\otimes\Res_{G_P}^G(\check{V})]\rangle[V]\\&&
= \frac{1}{n}\sum_{P\in X}\sum_{d=1}^{e_P-1} (e_P-d)
[\Ind_{G_P}^G(\cE(P)\otimes \chi_P^{-d})] \\ && = \frac{1}{n}
\sum_{P\in X} \sum_{d=1}^{e_P-1} d [\Ind_{G_P}^G(\cE(P)\otimes
\chi_P^d)].
\end{eqnarray*}
Thus, Corollary~1.3 follows from Theorem~1.1. \hfill $\Box$

The following corollary is the main result of the paper \cite{Na}
by Nakajima. Part (a) of it has also been proved by Kani (see
Theorem 2 in \cite{Ka}). In addition to Theorem~1.1 we use the
facts that the Euler characteristic $\chi(G,X,\cE)$ is an element
of the Grothendieck group $K_0(k[G])$ of projective $k[G]$-modules
(see Theorem~1 in \cite{Na} or Remark~1.5(a) below) and that the
Cartan homomorphism $K_0(k[G]) \ra K_0(G,k)$ is injective. The
corollary expresses the Euler characteristic $\chi(G,X,\cE)$ as an
{\em integral} linear combination of certain projective
$k[G]$-modules. Our proof shortens the somewhat lengthy
calculations in \cite{Na}.

{\bf Corollary 1.4} (Theorem 2 in \cite{Na}). \\ {\em (a) There is
a projective $k[G]$-module $N_{G,X}$ (which is unique up to
isomorphism) such that
\[\oplusm^n N_{G,X} \cong \oplusm_{P\in X} \oplusm_{d=0}^{e_P-1} \oplusm^d
\Ind_{G_P}^G(\chi_P^d).\] (b) For any $P \in X$, let
$l_{P,1},\ldots,l_{P,r} \in \{0,\ldots,e_P-1\}$ be given by the
equation
\[[\cE(P)] = \sum_{i=1}^r[\chi_P^{-l_{P,i}}] \quad \textsl{in}
\quad K_0(G_P,k);\] furthermore, for any $Q \in Y$, we fix a point
${\tilde{Q}} \in X$ with $\pi({\tilde{Q}}) = Q$. Then we have:
\[\chi(G,X,\cE) \equiv -r[N_{G,X}] + \sum_{Q\in Y} \sum_{i=1}^r
\sum_{d=1}^{l_{{\tilde{Q}},i}}
[\Ind_{G_{\tilde{Q}}}^G(\chi_{\tilde{Q}}^{-d})] \textsl{ in }
K_0(k[G]) \textrm{ mod } \ZZ[k[G]].\]}

{\bf Proof}. \\ (a) Applying Theorem 1.1 to the sheaf $\cE =
\cO_X$ with trivial $G$-action, we obtain the following equality
in $K_0(G,k)$ and hence in $K_0(k[G])$:
\[n \cdot \chi(G,X,\cO_X) = n(1-g_Y)[k[G]] - \sum_{P\in X}
\sum_{d=0}^{e_P-1} d[\Ind_{G_P}^G(\chi_P^d)].\] This equality
shows that the class of the projective $k[G]$-module
$\oplus_{d=0}^{e_P-1} \oplus^d \Ind_{G_P}^G(\chi_P^d)$ is
divisible by $n$ in $K_0(k[G])$. Writing the quotient as a linear
combination of classes of indecomposable projective modules we see
that the quotient is in fact the class of a projective
$k[G]$-module, say $N_{G,X}$. This immediately implies part~(a).\\
(b) We first prove the following congruence:
\[\deg(\cE) \equiv \sum_{P\in X}\sum_{i=1}^r l_{P,i} \textrm{ mod } n.\]
For this, we may obviously assume that $r=1$. We write $l_P$ for
$l_{P,1}$. Let $\cal K$ denote the sheaf of meromorphic functions
on $X$, i.e., the constant sheaf associated with the function
field $K(X)$ of $X$. Then $\cE$ is a $G$-subsheaf of the constant
sheaf $\cE \otimes \cK$. But $\cE \otimes \cK$ is isomorphic to
$\cK$ as a $G$-sheaf since the twisted group ring $K(X)G$ is
isomorphic to the ring $M_n(K(Y))$ of $n \times n$-matrices over
the function field $K(Y)$ by Galois theory and since there is (up
to isomorphism) only one module of $K(Y)$-dimension $n$ over
$M_n(K(Y))$. So, we may assume that $\cE = \cO_X(D)$ for some
equivariant Weil divisor $D = \sum_{P\in X} n_P [P]$. Now it is
easy to see that $l_P \equiv n_p$ mod $e_P$ for all $P$. So, for
any $Q \in Y$ we have: $\sum_{P\in \pi^{-1}(Q)} l_P =
\frac{n}{e_P} l_P \equiv \frac{n}{e_P} n_P = \sum_{P\in
\pi^{-1}(Q)} n_P$ mod $n$. Thus, the congruence above is proved.\\
Hence, by Theorem 1.1, we have in $K_0(G,k)/n\ZZ[k[G]]$:
\begin{eqnarray*}
\lefteqn{n \cdot \chi(G,X,\cE)=}\\&& = \sum_{P\in X} \sum_{i=1}^r
l_{P,i} [k[G]] - \sum_{P\in X} \sum_{i=1}^r \sum_{d=1}^{e_P-1}
d[\Ind_{G_P}^G(\chi_P^{d-l_{P,i}})] \\&& = \sum_{P\in X}
\sum_{i=1}^r \Ind_{G_P}^G \left(\sum_{d=0}^{e_P-1}
l_{P,i}[\chi_P^d] - \sum_{d=0}^{e_P-1}
d[\chi_P^{d-l_{P,i}}]\right) \\ && = \sum_{P\in X} \sum_{i=1}^r
\Ind_{G_P}^G\left(\sum_{d=0}^{e_P-1}
\alpha_{P,i,d}[\chi_P^d]\right)
\end{eqnarray*}
where $\alpha_{P,i,d} := \left\{
\begin{array}{ll}
-d & \textrm{for } d=0,\ldots, e_P-l_{P,i}-1 \\ -d+e_P &
\textrm{for } d= e_P-l_{P,i}, \ldots, e_P-1.
\end{array}\right.$\\
On the other hand, we have:
\begin{eqnarray*}
\lefteqn{-n \cdot r \cdot [N_{G,X}]+ n \sum_{Q\in Y} \sum_{i=1}^r
\sum_{d=1}^{l_{{\tilde{Q}},i}}[\Ind_{G_{\tilde{Q}}}^G(\chi_{\tilde{Q}}^{-d})]=}\\&&
= -r \sum_{P\in X} \sum_{d=0}^{e_P-1} d[\Ind_{G_P}^G (\chi_P^d)] +
\sum_{P\in X} e_P\sum_{i=1}^r
\sum_{d=1}^{l_{P,i}}[\Ind_{G_P}^G(\chi_P^{-d})]\\&&= \sum_{P\in X}
\sum_{i=1}^r \Ind_{G_P}^G\left(-\sum_{d=0}^{e_P-1} d[\chi_P^d] +
\sum_{d=1}^{l_{P,i}} e_P[\chi_P^{e_P-d}]\right)
\\&& = \sum_{P\in
X} \sum_{i=1}^r \Ind_{G_P}^G\left(\sum_{d=0}^{e_P-1}
\alpha_{P,i,d}[\chi_P^d]\right).
\end{eqnarray*}
Thus, Corollary 1.4 is proved. \hfill $\Box$

{\bf Remark 1.5}. \\ (a) In order to prove that $\chi(G,X,\cE)$ is
in $K_0(k[G])$ for all locally free $G$-modules $\cE$, it suffices
to show that the element $\frac{1}{n} \sum_{P\in X}
\sum_{d=0}^{e_p-1} d [\Ind_{G_P}^G(\chi_P^d)]$ is in $K_0(k[G])$
which in turn follows from the fact that $\chi(G,X,\cL)$ is in
$K_0(k[G])$ for {\em one} invertible $G$-module $\cL$ on~$X$.
(Apply twice the formula in part (b) of Corollary 1.4). If
$\deg(\cL) > 2g_X-2$, a nice and short proof of this fact using
equivariant cohomology can be found in Borne's thesis (see
Corollaire~3.14 on
p.~61 in \cite{Bo}). \\
(b) The equation
\[\chi(G,X,\cO_X) = (1-g_Y)[k[G]] - [N_{G,X}]\]
occurring in the proof of part (a) may be considered as an
equivariant version of the classical Hurwitz formula, see
Th\'eor\`eme 3.16 on p.\ 62 in Borne's thesis \cite{Bo}. He gives
a proof of this formula and of part~(a) of Corollary~1.4 which
does not use the work of Nakajima or Kani either and whose main
ingredient is the coherent Lefschetz fixed point formula as well.

The following example illustrates part (a) of Corollary 1.4; it
has been proved in \cite{Bo} directly using Hilbert 90 (see
Proposition 3.7 on p.\ 56 in \cite{Bo}).

{\bf Example 1.6}. We assume that $e_P = n$ for all $P \in X_{\rm
ram} := \{P\in X: e_P \not= 1\}$. Let $\chi : G \ra k^\times$ be a
character. We write $\chi = \chi_P^{a_P}$ for some $a_P \in \NN$
(for all $P\in X_{\rm ram}$). Then we have:
\[\sum_{P\in X_{\rm ram}} a_P \equiv 0 \textrm{ mod } n.\]

The following corollary is a main result of the paper \cite{Ka} by
Kani; it generalizes the classical Chevalley-Weil formula.

{\bf Corollary 1.7} (Corollary of Theorem 2 in \cite{Ka}). {\em
Let $\Omega$ denote the sheaf of holomorphic differentials on $X$.
Then we have in $K_0(G,k)$:
\[[H^0(X,\Omega)] = [k] + (g_Y-1)[k[G]] + [\check{N}_{G,X}].\] }

{\bf Proof}. It is well-known that $\deg(\Omega) = 2g_X-2$ and
that $\Omega(P)$ is $k[G_P]$-isomorphic to ${\mathfrak
m}_P/{\mathfrak m}_P^2$ for all $P\in X$. Hence, by Theorem 1.1,
we have in $K_0(G,k)$:
\begin{eqnarray*}
\lefteqn{n \cdot \chi(G,X,\Omega)=}\\&& = \left(n(1-g_Y) +
(2g_X-2) \right)[k[G]] - \sum_{P\in X} \sum_{d=0}^{e_P-1}
d[\Ind_{G_P}^G(\chi_P^{d+1})]\\&& = \left(n(1-g_Y) + n(2g_Y-2)
\right)[k[G]] - \sum_{P\in X} \sum_{d=0}^{e_P-1}
(d+1-e_P)[\Ind_{G_P}^G(\chi_P^{d+1})]\\&& = n(g_Y-1)[k[G]] +
n[\check{N}_{G,X}].
\end{eqnarray*}
Since $H^1(X,\Omega)$ is isomorphic to the trivial representation
$k$, this proves Corollary 1.7. \hfill $\Box$

The following corollary generalizes a recently published result of
Ksir (see \cite{Ks}). While we derive it from the previous
corollary, her proof is much more elementary. To be more precise:
While our proof is based on the Lefschetz fixed point formula (see
the proof of Theorem~1.1), her proof uses only the Riemann-Roch
and Hurwitz theorem and some elementary character theory. However,
her proof seems to work only in the case that not only the
representation $V$, but all (irreducible) $k$-representations of
$G$ are rationally valued. Here, we call a $k[G]$-module {\em
rationally valued} if its (Brauer) character takes only rational
values. For each point $Q \in Y$ , we fix a point $\tilde{Q} \in
X$ in the fibre $\pi^{-1}(Q)$.

{\bf Corollary~1.8}. {\em We assume that ${\rm char}(k)$ does not
divide $n$. Let $V$ be a non-trivial rationally valued irreducible
$k$-representation of $G$. Then the multiplicity of $V$ in the
$k$-representation $H^0(X,\Omega_X)$ is equal to
\[\dim(V)(g_Y-1) + \frac{1}{2}\sum_{Q\in Y} \left(\dim(V) -
\dim(V^{G_{\tilde{Q}}})\right).\] }

{\em Proof.} This follows from Corollary~1.7, the following
Proposition~1.9 and the well-known fact that the multiplicity of
$V$ in the regular representation $k[G]$ is equal to ${\rm
dim}(V)$.

{\bf Proposition~1.9}. {\em We assume that ${\rm char}(k)$ does
not divide $n$. Let $V$ be a rationally valued irreducible
$k$-representation of $G$. Then the multiplicity of $V$ in
$N_{G,X}$ and in its dual $\check{N}_{G,X}$ is equal to
\[\frac{1}{2} \sum_{Q\in Y} \left(\dim(V) -
\dim(V^{G_{\tilde{Q}}})\right).\]}

{\em Proof.} As in the proof of Corollary~1.3 we write $\langle \,
, \, \rangle$ for the usual character pairing. Then we have:
\begin{eqnarray*}
\lefteqn{\langle [V], [N_{G,X}]\rangle = \frac{1}{n} \sum_{P\in X}
\sum_{d=0}^{e_P-1} d \, \langle[V],[\Ind_{G_P}^G(\chi_P^d)]\rangle
}\\
&=& \frac{1}{n} \sum_{P\in X} \sum_{d=0}^{e_P-1} d \, \langle
[\Res^G_{G_P}(V)], [\chi_P^d]\rangle \quad \textrm{(by Frobenius
reciprocity)} \\
&=& \frac{1}{n} \sum_{P\in X} \frac{e_P}{2} \left(\dim(V) -
\dim(V^{G_P})\right) \quad \textrm{(by Lemma~1.10 below)}\\
&=& \frac{1}{2} \sum_{Q \in Y} \left( \dim(V) -
\dim(V^{G_{\tilde{Q}}})\right)
\end{eqnarray*}
since over any point $Q \in Y$ there are precisely
$\frac{n}{e_{\tilde{Q}}}$ points in the fibre $\pi^{-1}(Q)$ and
since $G_{\tilde{Q'}}$ is conjugate to $G_{\tilde{Q}}$ for any
other point $\tilde{Q'}$ in $\pi^{-1}(Q)$. This proves the
Proposition for $N_{G,X}$. The same argument applies to
$\check{N}_{G,X}$. \hfill $\Box$

{\bf Lemma~1.10.} {\em Let $C$ be a cyclic group of order $c$
coprime to $p$, let $V$ be a rationally valued $k[C]$-module, and
let $\chi:C \ra k^\times$ be a primitive character of $C$. Then we
have:
\[\sum_{d=0}^{c-1} d \, \langle [\chi^d],[V] \rangle = \frac{c}{2} \,
\left(\dim(V) - \dim(V^C)\right) \quad \textrm{in} \quad \QQ.\]}

{\em Proof.} It obviously suffices to consider the case $k = \CC$.
Since $V$ is rationally valued and $C$ is abelian, the class $[V]$
of $V$ in $K_0(\CC[C])$ belongs to the image of the canonical
homomorphism $K_0(\QQ[C]) \ra K_0(\CC[C])$, see the Corollary of
Proposition~35, \S 12.2, on p.~93 in Serre's book \cite{Se2}. By
Exercise~13.1 on pp.~104-105 in \cite{Se2}, the permutation
representations $\Ind_H^C(1_H)$, $H$ a subgroup of $C$, form a
$\ZZ$-basis of $K_0(\QQ[G])$. Since both sides of the formula in
the Lemma are additive in $V$, it therefore suffices to prove the
Lemma in the case $V=\Ind_H^C(1_H)$ where $H$ is any subgroup of
$C$. Let $h$ denote the order of $H$. Then we obviously have:
\[\langle[\chi^d],[\Ind_H^C(1_H)]\rangle = \langle
[\Res_H^C(\chi^d)],[1_H] \rangle = \left\{
\begin{array}{ll} 1, & \textrm{ if } h|d \\ 0, &\textrm{ else.}
\end{array}\right.\]
Thus we obtain:
\begin{eqnarray*}
\lefteqn{ \sum_{d=0}^{c-1} d \, \langle [\chi^d],
[\Ind_H^C(1_H)] \rangle }\\
&=&  h + 2h + \ldots + \left(\frac{c}{h}-1\right)h\\
&=& h \left(1 + 2 + \ldots + \left(\frac{c}{h} -1\right)\right)\\
&=& h \frac{(\frac{c}{h}-1) \frac{c}{h}}{2} =
\frac{c}{2}\left(\frac{c}{h} -1\right) \\
&=& \frac{c}{2} \left(\dim(\Ind_H^C(1_H)) -
\dim\left((\Ind_H^C(1_H))^C\right)\right),
\end{eqnarray*}
as was to be shown. \hfill $\Box$

Similarly to the deduction of Corollary~1.8 from Corollary~1.7 we
deduce the following corollary from Corollary~1.4. An alternative
approach to the following corollary based on Ksir's paper
\cite{Ks} and Borne's thesis \cite{Bo} ca be found in the recent
preprint \cite{JK} by Ksir and Joyner.

{\bf Corollary~1.11}. {\em We assume that ${\rm char}(k)$ does not
divide $n$. Let $D= \sum_{P\in X} n_P [P]$ be a $G$-equivariant
divisor on $X$ and let $V$ be a rationally valued irreducible
$k$-representation of $G$. Then the multiplicity of $V$ in the
Euler characteristic $\chi(G,X,\cO_X(D))$ is equal to
\begin{eqnarray*}\lefteqn{\dim(V)\left(1-g_Y\right) + }\\
&& \sum_{Q \in Y} \left(\dim(V)m_{\tilde{Q}} - \frac{1}{2}
 \left(\dim(V) - \dim(V^{G_{\tilde{Q}}})\right) +
 \sum_{d=1}^{l_{\tilde{Q}}} \left\langle
[\chi_{\tilde{Q}}^{-d}], [\Res_{G_{\tilde{Q}}}^G(V)] \right\rangle
\right)
\end{eqnarray*} where $l_P\in \{0, \ldots, e_P-1\}$ and
$m_P \in \ZZ$ are given by $n_P = l_P + m_P e_P$ for any $P \in
X$.}

{\em Proof.} By Corollary~(1.4)(b) we have the congruence
\[\chi(G,X,\cO_X(D)) \equiv - [N_{G,X}] + \sum_{Q \in Y}
\sum_{d=1}^{l_{\tilde{Q}}} [\Ind_{G_{\tilde{Q}}}^G(\chi_P^{-d})]
\] in $K_0(k[G])$ mod $\ZZ[k[G]]$. By the Riemann-Roch theorem and
the Riemann-Hurwitz formula this congruence becomes an equality in
$K_0(k[G])$ after adding the term $\left(1-g_Y + \sum_{Q \in
Y}m_{\tilde{Q}}\right)[k[G]]$ on the right hand side. Now using
the Proposition~1.9 and Frobenius reciprocity we obtain Theorem~2.
\hfill $\Box$

{\em Remark}. Note that the multiplicity $\langle
[\chi_{\tilde{Q}}^{-d}], [\Res_{G_{\tilde{Q}}}(V)]\rangle$ of the
character $\chi_{\tilde{Q}}^{-d}$ of the cyclic group
$G_{\tilde{Q}}$ in the restricted representation
$\Res_{G_{\tilde{Q}}}^G(V)$ can easily be computed. For instance,
if $e_{\tilde{Q}}$ is a power of a prime, the following lemma can
be applied.

{\bf Lemma~1.12.} {\em Let $C$ be a cyclic group of prime power
order $l^r$, let $V$ be a rationally valued $k$-representation of
$C$, let $\chi: C \ra k^\times$ be a primitive character and let
$d=l^s m \in \ZZ$ with $s \in \{0, \ldots, r-1\}$ and $m \in \ZZ$
coprime to $l$. Then the multiplicity $\langle [\chi^d], [V]
\rangle$ of $\chi^d$ in $V$ is equal to
\[ \frac{\dim(V^H) - \dim(V^{H'})}{l^{r-s-1}(l-1)}\]
where $H$ and $H'$ are the (unique) subgroups of $C$ of order
$l^s$ and $l^{s+1}$, respectively.}

{\em Proof.} Easy. \hfill $\Box$

\bigskip

\section*{\S 2 The equivariant Grothendieck-Ogg-Shafarevich
formula}

The goal of this section is to prove an equivariant
Grothendieck-Ogg-Shafarevich formula.

Let $X$ be a connected smooth projective curve over an
algebraically closed field $k$ and let $G$ be a finite subgroup of
$\Aut(X/k)$ of order $n$. We assume in this section that ${\rm
char}(k)$ does not divide $n$. Let $\pi:X\ra Y$, $g_X$, $g_Y$,
$G_P$ and $e_P$ be defined as in \S 1. Furthermore, we denote the
generic point of $X$ by $\eta$.

Let $l\not= {\rm char}(k)$ be a prime and let $\cF$ be a
constructible $\FF_l$-sheaf on $X_{\rm \acute{e}t}$ with
$G$-action, i.e., we have isomorphisms $g^*(\cF) \ra \cF$, $g\in
G$, which satisfy the usual composition rules. Then the \'etale
cohomology groups $H^i_{\rm \acute{e}t}(X,\cF)$, $i\ge 0$, are
$\FF_l$-representations of $G$. Let
\[\chi_{\rm \acute{e}t}(G,X,\cF) :=  [H_{\rm \acute{e}t}^0(X,\cF)] - [H_{\rm
\acute{e}t}^1(X,\cF)]+ [H_{\rm \acute{e}t}^2(X,\cF)] \in
K_0(G,\FF_l)\] denote the equivariant Euler characteristic; here,
$K_0(G,\FF_l)$ is the Grothendieck group of
$\FF_l$-representations of $G$ (of finite dimension). Furthermore,
let
\[\alpha(\cF) := \sum_{P\in X} \alpha_P(\cF) \in \ZZ\]
denote the sum of the wild conductors of $\cF$ (see p.\ 188 in
\cite{Mi}) and let $I_P:= {\rm ker}(\FF_l[G_P] \ra \FF_l)$ denote
the augmentation representation of $G_P$ (for $P \in X$).

The following theorem may be viewed as an analogue of Theorem 1.1;
it computes the equivariant Euler characteristic $\chi_{\rm
\acute{e}t}(G,X,\cF)$.

{\bf Theorem 2.1} (Equivariant Grothendieck-Ogg-Shafarevich
formula). {\em We have in $K_0(G,\FF_l)_\QQ$:
\begin{eqnarray*}
\lefteqn{\chi_{\rm \acute{e}t}(G,X,\cF) = }\\& = &
\left((2-2g_Y)\dim_{\FF_l}(\cF_{\bar{\eta}}) - \frac{1}{n}
\alpha(\cF) + \sum_{P\in X} \frac{e_P}{n}
\left(\dim_{\FF_l}(\cF_P) -
\dim_{\FF_l}(\cF_{\bar{\eta}})\right)\right)[\FF_l[G]] \\ &&
-\sum_{P\in X} \frac{e_P}{n} [\Ind_{G_P}^G(\cF_P \otimes I_P)].
\end{eqnarray*}}

{\bf Remark 2.2}. Let $\dim_{\FF_l}(\cF_P) =
\dim_{\FF_l}(\cF_{\bar{\eta}})$ for all $P \in X$ with $e_P \not=
1$. Then the formula above has the following shape:
\[\chi_{\rm \acute{e}t}(G,X,\cF) = \left((2-2g_Y)
\dim_{\FF_l}(\cF_{\bar{\eta}}) - \frac{1}{n} c(\cF)\right)
[\FF_l[G]] - \sum_{P\in X} \frac{e_P}{n}
[\Ind_{G_P}^G(\cF_P\otimes I_P)];\] here, $c(\cF) = \sum_{P\in X}
c_P(\cF)$ is the sum of the conductors of $\cF$ (see p.\ 188 in
\cite{Mi}). This formula becomes particularly simple in the
following two extreme cases.\\(a) Let $\pi: X \ra Y$ be \'etale.
Then we have:
\[\chi_{\rm \acute{e}t}(G,X,\cF) = \left( (2-2g_Y)
\dim_{\FF_l}(\cF_{\bar{\eta}}) - \frac{1}{n} c(\cF)\right)
[\FF_l[G]].\]  If $G$ is the trivial group, this is the classical
Grothendieck-Ogg-Shafarevich formula (see Theorem 2.12 on p.\ 190
in \cite{Mi}). In particular, we obtain the following formula for
an arbitrary group $G$:
\[\chi_{\rm \acute{e}t}(G,X,\cF) = \frac{1}{n} \cdot \chi_{\rm \acute{e}t}(X,\cF)
\cdot [\FF_l[G]]\] (note that $n(2-2g_Y) = (2-2g_X)$ by the
Hurwitz formula). This formula remains valid, if we drop the
assumption that ${\rm char(k)}$ does not divide $n$ (see the proof
below). In particular we obtain that the (non-equivariant) Euler
characteristic $\chi_{\rm \acute{e}t}(X,\cF)$ is divisible by $n$.
Finally, the latter formula may be viewed as an analogue of
Theorem 2.4 in \cite{EL}.\\ (b) Let $\cF$ be the constant sheaf
$\FF_l$ with trivial $G$-action. Then we obtain the following
formula:
\[\chi_{{\rm \acute{e}t}}(G,X,\cF) = (2-2g_Y)[\FF_l[G]] - \sum_{P\in
X} \frac{e_P}{n}[\Ind_{G_P}^G(I_P)].\] This is the $\FF_l$-version
of the formula in Remark 2.9 on p.\ 187 in \cite{Mi}. (Note that
the Artin character is the character of the augmentation
representation since ${\rm char}(k)$ does not divide $n$.) It can
be derived from the $\QQ_l$-version by applying the decomposition
homomorphism as in lines 6 through 9 on p.\ 191 in \cite{Mi}.

{\bf Proof} (of Theorem 2.1). As in the proof of Theorem~1.1 we
will show that the (Brauer) character values of both sides
coincide for all $\sigma \in G$. So let $\sigma \in G$. Then we
have:
\begin{eqnarray*}
\lefteqn{\frac{1}{n} \sum_{P\in X} e_P \cdot
\Trace(\sigma|\Ind_{G_P}^G(\cF_P \otimes I_P)) = }\\ && =
\frac{1}{n}\sum_{P\in X} \sum_{\tau \in G, \, \tau^{-1}\sigma\tau
\in G_P} \Trace(\tau^{-1}\sigma\tau|\cF_P\otimes I_P)\\ &&
\hspace*{5cm} \textrm{(by Lemma~(21.28) on p.~509 in
\cite{CR2})} \\
&& = \frac{1}{n} \sum_{P\in X} \sum_{\tau\in G,\, \tau(P) \in
X^{\sigma}} \Trace(\sigma|\cF_{\tau(P)} \otimes I_{\tau(P)})\\&& =
\sum_{P\in X^{\sigma}} \Trace(\sigma|\cF_P \otimes I_P) \\ && =
\sum_{P\in X^{\sigma}} \left(\Trace(\sigma|
\cF_P\otimes\FF_l[G_P]) - \Trace(\sigma|\cF_P)\right)\\&& =
\sum_{P\in X^\sigma} \left(\dim_{\FF_l}(\cF_P) \cdot
\Trace(\sigma|\FF_l[G_P]) - \Trace(\sigma|\cF_P)\right)
\textrm{(by Frobenius reciprocity)} \\ && = \left\{
\begin{array}{ll}
- \summ_{P\in X^\sigma} \Trace(\sigma|\cF_P) & \textrm{for }
\sigma \not= \id \\ \summ_{P\in X} \dim_{\FF_l}(\cF_P)\cdot(e_P-1)
& \textrm{for } \sigma = \id .
\end{array}\right.
\end{eqnarray*}
Hence, for $\sigma \not= \id$, the character value of the right
hand side of the formula in Theorem~2.1 at the place $\sigma$
equals $\sum_{P\in X^\sigma} \Trace(\sigma|\cF_P)$. By the
Lefschetz fixed point formula (see Theorem 2 in \cite{Do2} or
\cite{Ve}; here, we use that ${\rm char}(k)$ does not divide $n$),
this equals the character value of the left hand side.\\ For
$\sigma = \id$, the character value of the right hand side at the
place $\sigma$ is
\begin{eqnarray*}
\lefteqn{n(2-2g_Y)\dim_{\FF_l}(\cF_{\bar{\eta}})- \alpha(\cF) +
\sum_{P \in X} \dim_{\FF_l}(\cF_P) - \sum_{P\in X} e_P \cdot
\dim_{\FF_l}(\cF_{\bar{\eta}}) = }\\ && = (2-2g_X)
\dim_{\FF_l}(\cF_{\bar{\eta}}) - \left(\alpha(\cF) + \sum_{P\in X}
(\dim_{\FF_l}(\cF_{\bar{\eta}}) - \dim_{\FF_l}(\cF_P))\right)
\end{eqnarray*}
by the Hurwitz formula. By the classical
Grothendieck-Ogg-Shafarevich formula (see Theorem 2.12 on p.\ 190
in \cite{Mi}), this equals the character value of the left hand
side at the place $\sigma = \id$. \\ Thus, the proof of Theorem
2.1 is complete. \hfill $\Box$

The following corollary may be viewed as the analogue of Corollary
1.3; it computes the multiplicity of any irreducible
$\FF_l$-representation $V$ of $G$ in the equivariant Euler
characteristic $\chi_{\rm \acute{e}t}(G,X,\cF)$. We write ${\rm
Irr}(G,\FF_l)$ for the set of isomorphism classes of irreducible
$\FF_l$-representations and set $s_V:=
\dim_{\FF_l}(\End_{\FF_l[G]}(V))$ for $V \in {\rm Irr}(G,\FF_l)$.
For $Q\in Y$ and $V \in {\rm Irr}(G,\FF_l)$, let $m_{Q,\cF,V}$
denote the multiplicity of the trivial representation $\FF_l$ in
$\cF_P \otimes I_P \otimes \Res_{G_P}^G(V)$ where $P \in
\pi^{-1}(Q)$.

{\bf Corollary 2.3}. {\em We assume that $(l, \ord(G)) = 1$. Then
we have in $K_0(G,\FF_l)_\QQ$:
\begin{eqnarray*}
\lefteqn{\chi_{\rm \acute{e}t}(G,X,\cF) = }\\ & = &
\left((2-2g_Y)\dim_{\FF_l}(\cF_{\bar{\eta}}) - \frac{1}{n}
\alpha(\cF) + \sum_{P \in X} \frac{e_P}{n}(\dim_{\FF_l}(\cF_P) -
\dim_{\FF_l}(\cF_{\bar{\eta}}))\right) [\FF_l[G]] \\ && -
\sum_{V\in {\rm Irr}(G,\FF_l)} \frac{1}{s_V}\sum_{Q\in Y}
m_{Q,\cF,\check{V}}[V].
\end{eqnarray*}}

{\bf Proof}. Let $\langle, , \,\rangle: K_0(G,\FF_l) \times
K_0(G,\FF_l) \ra \ZZ$ denote the symmetric bilinear form given by
$\langle[V],[W]\rangle := \dim_{\FF_l}(\Hom_{\FF_l[G]}(V,W))$ for
any
$\FF_l$-representations $V$, $W$ of $G$. Then we obviously have:\\
(a) $\langle[V],[\check{W} \otimes X]\rangle = \langle[V\otimes
W], [X]\rangle$ for any $\FF_l$-representations $V$, $W$, $X$ of $G$. \\
(b) $\sum_{V \in {\rm Irr}(G,\FF_l)}
\frac{1}{\langle[V],[V]\rangle} \langle[V],x\rangle [V] = x$ for
all $x \in K_0(G,\FF_l)$. \\ (c) $\langle x,\Res_H^G(y)\rangle =
\langle\Ind_H^G(x), y\rangle$ for any subgroup
$H$ of $G$, $x \in K_0(H,\FF_l)$ and $y \in K_0(G,\FF_l)$. \\
Hence we have for all $P\in X$:
\begin{eqnarray*}
\lefteqn{\sum_{V \in {\rm Irr}(G,\FF_l)}
\frac{1}{\langle[V],[V]\rangle} \langle[\FF_l],[\cF_P \otimes I_P
\otimes \Res_{G_P}^G(\check{V})]\rangle[V]=}
\\&& = \sum_{V \in {\rm Irr}(G,\FF_l)} \frac{1}{\langle[V],[V]\rangle}
\langle[V],[\Ind_{G_P}^G(\cF_P \otimes I_P)]\rangle[V] \\ && =
[\Ind_{G_P}^G(\cF_P\otimes I_P)] \textrm{ in } K_0(G,\FF_l).
\end{eqnarray*}
Hence we have:
\begin{eqnarray*}
\lefteqn{\sum_{V\in {\rm Irr}(G,\FF_l)} \frac{1}{s_V} \sum_{Q\in
Y} m_{Q,\cF,\check{V}}[V] =} \\ && = \sum_{P\in X} \frac{e_P}{n}
\sum_{V \in {\rm Irr}(G,\FF_l)} \frac{1}{\langle[V],[V]\rangle}
\langle[\FF_l],[\cF_P \otimes I_P \otimes \Res_{G_P}^G(\check{V})]\rangle[V]\\
&& =  \sum_{P\in X}\frac{e_p}{n}[\Ind_{G_P}^G (\cF_P \otimes I_P)]
\textrm{ in } K_0(G,\FF_l).
\end{eqnarray*}
Thus, Corollary~2.3 follows from Theorem~2.1. \hfill $\Box$

The analogue of the element $[N_{G,X}]$ occurring in Corollary
1.4(a) is $$\sum_{P\in X} \frac{e_P}{n} [\Ind_{G_P}^G(I_P)].$$ It
is obviously an element of $K_0(G,\FF_l)$. More generally,
$\sum_{P\in X} \frac{e_P}{n} [\Ind_{G_P}^G(\cF_P\otimes I_P)]$ is
an element of $K_0(G,\FF_l)$, since $\Ind_{G_P}^G(\cF_P\otimes
I_P)$ is isomorphic to $\Ind_{G_{P'}}^G(\cF_{P'}\otimes I_{P'})$
for any $P, P' \in X$ with $\pi(P) = \pi(P')$. Furthermore,
$\sum_{P\in X} \frac{e_P}{n}\left(\dim_{\FF_l}(\cF_P) -
\dim_{\FF_l}(\cF_{\bar{\eta}})\right)$ is an integer. Thus,
Theorem 2.1 implies the following corollary.

{\bf Corollary 2.4}. {\em The sum $\alpha(\cF)$ of the wild
conductors $\alpha_P(\cF)$, $P\in X$, is divisible by~$n$.}

In particular, Theorem 2.1 expresses $\chi_{\rm
\acute{e}t}(G,X,\cF)$ as an {\em integral} linear combination of
$\FF_l$-representations; thus, the analogue of Corollary 1.4 is
already built into Theorem~2.1.

\bigskip

School of Mathematics\\
University of Southampton\\
Southampton SO17 1BJ\\
United Kingdom\\
{\em e-mail:} bk@maths.soton.ac.uk.


\bigskip

{\footnotesize
\begin{thebibliography}{BFQ}

\bibitem[BFQ]{BFQ} {\sc P.\ Baum, W.\ Fulton} and {\sc G.\ Quart},
Lefschetz-Riemann-Roch for singular varieties, {\em Acta Math.}\
{\bf 143} (1979), 193-211.
\bibitem[Bo]{Bo} {\sc N.\ Borne}, Une formule de Riemann-Roch
\'equivariante pour les courbes, {\em Th\`ese} (Universit\'e
Bordeaux I, 2000).
\bibitem[CR2]{CR2} {\sc C.~W.~Curtis} and {\sc I.~Reiner}, ``Methods
of representation theory, with applications to finite groups and
orders'', Vol.~I, {\em Pure and Applied Mathematics}, John Wiley
\& Sons, New York (1981).
\bibitem[Do1]{Do1} {\sc P.\ Donovan}, The Lefschetz-Riemann-Roch
formula, {\em Bull.\ Soc.\ Math.\ France}~{\bf 97} (1969),
257-273.
\bibitem[Do2]{Do2} {\sc P.\ Donovan}, The fixed point formula for
equivariant \'etale sheaves, {\em J.\ Algebra}~{\bf 103} (1986),
427-436.
\bibitem[EL]{EL} {\sc G.\ Ellingsrud} and {\sc K.\ L{\o}nsted}, An
equivariant Lefschetz formula for finite reductive groups, {\em
Math.\ Ann.}\ {\bf 251} (1980), 253-261.
\bibitem[Ha]{Ha} {\sc R.\ Hartshorne}, ``Algebraic geometry",
{\em Graduate Texts in Mathematics}~{\bf 52}, Springer-Verlag, New
York (1977).
\bibitem[JK]{JK} {\sc A.~Ksir} and {\sc D.~Joyner},
Representations of finite groups on Riemann-Roch spaces, II, {\em
preprint} (2003), 16 pp.
\bibitem[Ka]{Ka} {\sc E.\ Kani}, The Galois-module structure of
the space of holomorphic differentials of a curve, {\em J.\ Reine
Angew.\ Math.}\ {\bf 367} (1986), 187-206.
\bibitem[K\"o1]{Ko} {\sc B.\ K\"ock}, The Lefschetz theorem in
higher equivariant $K$-theory, {\em Comm.\ Algebra} {\bf 19}
(1991), 3411-3422.
\bibitem[K\"o2]{Ko2} {\sc B.~K\"ock}, Galois structure of Zariski
cohomology for weakly ramified covers of curves, to appear in {\em
Amer.\ J.\ Math.}~{\bf 126}, 23 pp.
\bibitem[Ks]{Ks} {\sc A.~E.~Ksir}, Dimensions of Prym varieties,
{\em Int.\ J.\ Math.\ Sci.}~{\bf 26} (2001), 107-116.
\bibitem[Mi]{Mi} {\sc J.\ S.\ Milne}, ``\'Etale Cohomology", {\em
Princeton Mathematical Series} {\bf 33}, Princeton University
Press, Princeton (1980).
\bibitem[Na]{Na} {\sc S. Nakajima}, Galois module structure of
cohomology groups for tamely ramified coverings of algebraic
varieties, {\em J.\ Number Theory} {\bf 22} (1986), 115-123.
\bibitem[Se1]{Se1} {\sc J.-P.\ Serre}, ``Corps locaux'', {\em
Publications de l'Institut de Math\'ematique de l'Universit\'e de
Nancago}~{\bf VIII}, Hermann, Paris (1962).
\bibitem[Se2]{Se2} {\sc J.-P.~Serre}, ``Linear representations of
finite groups'', {\em Graduate Texts in Mathematics}~{\bf 42},
Springer-Verlag, New York (1977).
\bibitem[Ve]{Ve} {\sc J.\ L.\ Verdier}, The Lefschetz fixed point
formula in \'etale cohomology, {\em in} ``Proceedings of a
conference on local fields" (T.\ A.\ Springer, Ed.),
Springer-Verlag, Berlin (1967), 199-214.


\end{thebibliography} }
\end{document}